\documentclass[a4paper,11pt]{amsart}

\newfont{\cyr}{wncyr10 scaled 1100}

\usepackage[left=2.7cm,right=2.7cm,top=3.5cm,bottom=3cm]{geometry}
\usepackage{amsthm,amssymb,amsmath,amsfonts,mathrsfs,amscd}
\usepackage[latin1]{inputenc}
\usepackage[all]{xy}
\usepackage{latexsym}
\usepackage{longtable}

\theoremstyle{plain}
\newtheorem{theorem}{Theorem}[section]
\newtheorem{corollary}[theorem]{Corollary}
\newtheorem{lemma}[theorem]{Lemma}
\newtheorem{proposition}[theorem]{Proposition}

\theoremstyle{definition}
\newtheorem{definition}[theorem]{Definition}

\theoremstyle{remark}

\newtheorem{remark}[theorem]{Remark}

\newcommand{\Emb}{\operatorname{Emb}}

\newcommand{\g}{\gamma}

\newcommand{\Q}{\mathbb{Q}}
\newcommand{\Z}{\mathbb{Z}}
\newcommand{\F}{\mathbb{F}}

\newcommand{\C}{\mathbb{C}}

\newcommand{\PP}{\mathbb{P}}

\newcommand{\Gal}{\operatorname{Gal}}
\newcommand{\GL}{\operatorname{GL}}

\newcommand{\Div}{\operatorname{Div}}

\newfont{\gotip}{eufb10 at 12pt}

\newcommand{\cO}{{\mathcal O}}

\newcommand{\om}{{\omega}}

\newcommand{\SL}{{\mathrm {SL}}}

\newcommand{\R}{{\mathbb R}}
\newcommand{\M}{{\mathrm{M}}}

\newcommand{\PGL}{{\mathrm{PGL}}}

\newcommand{\W}{\mathbb W}

\newcommand{\q}{\underline{q}}

\def \mint {\times \hskip -1.1em \int}

\newcommand{\T}{\mathbb T}

\DeclareMathOperator{\Hom}{Hom} 

\newcommand{\X}{\mathbb X}
\newcommand{\Y}{\mathbb Y}
\newcommand{\U}{\mathbb U}

\newcommand{\res}{\mathrm{res}}
\newcommand{\D}{\mathbb D}

\newcommand{\longmono}{\mbox{$\lhook\joinrel\longrightarrow$}}

\newcommand{\smallmat}[4]{\bigl(\begin{smallmatrix}#1&#2\\#3&#4\end{smallmatrix}\bigr)}

\include{thebibliography}

\begin{document}

\title[Quaternionic Darmon points]{Quaternionic Darmon points on $p$-adic tori\\and abelian varieties}
\author{Matteo Longo and Stefano Vigni}

\begin{abstract}
We prove formulas for the $p$-adic logarithm of quaternionic Darmon points on $p$-adic tori and modular abelian varieties over $\Q$ having purely multiplicative reduction at $p$. These formulas are amenable to explicit computations and are the first to treat Stark--Heegner type points on higher-dimensional abelian varieties.
\end{abstract}

\address{Dipartimento di Matematica Pura e Applicata, Universit\`a di Padova, Via Trieste 63, 35121 Padova, Italy}
\email{mlongo@math.unipd.it}
\address{Departament de Matem\`atica Aplicada II, Universitat Polit\`ecnica de Catalunya, C. Jordi Girona 1-3, 08034 Barcelona, Spain}
\email{stefano.vigni@upc.edu}



\maketitle

\section{Introduction} \label{intro}

Stark--Heegner points on elliptic curves over $\Q$ were defined by Henri Darmon in \cite{Dar} as conjectural analogues over (abelian extensions of) real quadratic fields of classical Heegner points. Since then, generalizations of these points to higher-dimensional (modular) abelian varieties have not been systematically investigated.

Building on techniques developed in \cite{LRV1}, quaternionic Darmon points on $p$-adic tori and Jacobians of Shimura curves over $\Q$ were introduced in \cite{LRV2}. Recently, a rationality result for projections of Darmon points to elliptic curves (which extends to a quaternionic context the theorem obtained by Bertolini and Darmon in \cite{BD-annals}) has been proved in \cite{LV-darmon}. Let us briefly review the arithmetic setting in which the theory of \cite{LRV1} and \cite{LRV2} takes place.

Let $X_0$ (repectively, $X_1$) denote the (compact) Shimura curve of discriminant $D>1$ and level $M$ (respectively, $Mp$) for coprime integers $D,M$ and a prime number $p\nmid DM$ (here $D$ is a square-free product of an even number of primes). Let $H$ be the maximal torsion-free quotient of the cokernel of the degeneracy map from $H_1(X_0,\Z)^2$ to $H_1(X_1,\Z)$. It turns out that $H$ is a free abelian group of rank equal to twice the dimension of the $p$-new quotient $J_1^{\text{$p$-new}}$ of the Jacobian $J_1$ of $X_1$. Let $\mathbb H$ be a non-zero torsion-free quotient of $H$ and let $K$ be a real quadratic field in which all the primes dividing $M$ (respectively, $Dp$) split (respectively, are inert). The $p$-adic tori on which Darmon points are naturally defined are quotients of tori of the form $T_{\mathbb H}(K_p):=K^\times_p\otimes\mathbb H$ by suitable lattices $L_{\mathbb H}$ contained in $T_{\mathbb H}(\Q_p):=\Q_p^\times\otimes\mathbb H$ and built, as in \cite{LRV1}, via group (co)homology and $p$-adic integration. So far, explicit computations with Darmon points (of conductor $1$) had been performed only in \cite{LV-darmon} in the case where $\mathbb H$ is associated with an elliptic curve over $\Q$ of conductor $MDp$. 

While leading to rationality results for quaternionic Darmon points on elliptic curves over $\Q$ that provide evidence for the conjectures formulated in \cite{LRV2}, the one-dimensional setting studied in \cite{LV-darmon} might obscure one of the features of the theory developed in \cite{LRV1} and \cite{LRV2}, that is, the possibility of defining Darmon points of arbitrary conductor on higher-dimensional (modular) abelian varieties. It is perhaps worthwhile to remark that constructions of this kind lie outside the scope of Darmon's original theory, since Stark--Heegner points are defined directly on elliptic curves over $\Q$ (when Shimura curves reduce to classical modular curves, constructions analogous to ours could be made by using Dasgupta's $p$-adic uniformization of modular Jacobians, cf. \cite{Das}). 

The goal of the present paper is to obtain formulas for quaternionic Darmon points on $p$-adic tori and modular abelian varieties over $\Q$. To describe the main results in detail, we need some further notation. With $\mathbb H$ as before, let $\mathcal M_0$ denote the abelian group of $\mathbb H$-valued measures on $\PP^1(\Q_p)$ with total mass $0$. Following \cite[Section 4]{LRV1}, one can explicitly introduce a canonical cohomology class 
\[ \boldsymbol\mu_{\mathbb H}\in H^1(\Gamma,\mathcal M_0) \]
where $\Gamma\subset\SL_2(\Q_p)$ is Ihara's group defined in \eqref{Gamma}. The abelianization of $\Gamma$ is a finite group whose exponent will be denoted $t$. If $\cO_K$ is the ring of integers of $K$ and $d_K$ is its discriminant, let $\cO_c:=\Z+c\cO_K$ be the order of $K$ of conductor $c$ prime to $MDd_Kp$. As in \cite{LRV2}, multiplicative integration over $\PP^1(\Q_p)$ against a representative $\mu_{\mathbb H}$ of $\boldsymbol\mu_{\mathbb H}$ allows us to attach a Darmon point $\mathcal P_{\mathbb H,\psi}\in T_{\mathbb H}(K_p)/L_{\mathbb H}$ to every optimal embedding $\psi$ of $\cO_c$ into a fixed Eichler order $R_0$ of level $M$ in the quaternion algebra over $\Q$ of discriminant $D$. As conjectured in \cite[\S 3.2]{LRV2}, the points $\mathcal P_{\mathbb H,\psi}$ are expected to be rational over the narrow ring class field of $K$ of conductor $c$ and to satisfy a suitable Shimura reciprocity law under the action of the corresponding Galois group (see \cite{LV-darmon} for partial results in this direction when $c=1$).

Now let $\mathbb H_p:=\Z_p\otimes\mathbb H$, let $\mathbb X$ be the set of primitive vectors in $\mathbb Y:=\Z_p^2$ (i.e., those vectors which are not divisible by $p$) and write $\mathbb D$ for the group of $\mathbb H_p$-valued measures on $\mathbb Y$ that are supported on $\mathbb X$. If $\Gamma_0$ denotes the group of norm $1$ elements in $R_0$ and $\pi:\mathbb X\rightarrow\PP^1(\Q_p)$ is the map sending $(a,b)$ to $a/b$ then there exists $\boldsymbol{\tilde\mu}_{\mathbb H}\in H^1(\Gamma_0,\mathbb D)$ such that $\pi_*(\boldsymbol{\tilde\mu}_{\mathbb H})=\text{res}_{\Gamma_0}(\boldsymbol\mu_{\mathbb H})$. With self-explaining notation, we can also choose a $1$-cocycle $\tilde\mu_{\mathbb H}$ representing $\boldsymbol{\tilde\mu}_{\mathbb H}$ such that $\pi_*(\tilde\mu_{\mathbb H,\gamma})=\mu_{\mathbb H,\gamma}$ for all $\gamma\in\Gamma_0$. Let $\varepsilon_c$ denote a generator of the group of units of $\cO_c$ of norm $1$ such that $\varepsilon_c>1$ under a fixed embedding $K\hookrightarrow\R$ and set $\gamma_\psi:=\psi(\varepsilon_c)\in\Gamma_0$. Moreover, denote $z_\psi$ the fixed point of $\psi(K^\times)$ acting on $\PP^1(K_p)$ such that $\psi(\alpha)(z_\psi,1)=\alpha(z_\psi,1)$ for all $\alpha\in K$. Finally, fix an algebraic closure $\bar\Q_p$ of $\Q_p$ and let $m\geq1$ be the rank of $\mathbb H$. 

The first main result of this article, which we state below and is Theorem \ref{thm2.11} in the text, gives a formula for the value at $\mathcal P_{\mathbb H,\psi}$ of a continuous homomorphism
\[ \Psi:T_{\mathbb H}(\bar\Q_p)/L_{\mathbb H}\longrightarrow\bar\Q_p^m \]
that is $K_p^m$-valued on the subset $T_{\mathbb H}(K_p)/L_{\mathbb H}$ (the map $\Psi$ will also be viewed as defined on $\bar\Q_p^\times$ via diagonal embedding). 

\begin{theorem} \label{intro-thm1}
$\displaystyle{\Psi(\mathcal P_{\mathbb H,\psi})=-t\cdot\int_\X\Psi(x-z_\psi y)d{\tilde\mu}_{\mathbb H,\gamma_\psi}(x,y)}$. 
\end{theorem}  

By specializing Theorem \ref{intro-thm1} to a suitable choice of $\mathbb H$ and $\Psi$, we can also offer a formula for the images of Darmon points on abelian varieties over $\Q$ with purely multiplicative reduction at $p$, as now we explain. 

Let $A$ be a modular abelian variety over $\Q$ of dimension $d$ and conductor $N:=MDp$, which means that $A$ is associated with a (normalized) newform of weight $2$ and level $N$. It is known that $A$ has purely multiplicative reduction at $p$, i.e., the identity component of the Néron model of $A$ over $\Z_p$ is a torus over $\F_p$. Since $J_1^{\text{$p$-new}}$ is the maximal toric quotient of $J_1$ at $p$, it follows that $A$ is a quotient of $J_1^{\text{$p$-new}}$. Moreover, if $\epsilon\in\{\pm\}$ then the $\epsilon$-eigenspace for complex conjugation acting on $\mathbb H_A:=H_1(A(\C),\Z)$, which we denote $\mathbb H_A^\epsilon$, is a free quotient of $H$ of rank $d$. Set $T_A:=T_{\mathbb H_A^\epsilon}$, $L_A:=L_{\mathbb H_A^\epsilon}$ and $\tilde\mu_A:=\tilde\mu_{\mathbb H_A^\epsilon}$. As a consequence of the uniformization results of \cite{LRV1}, there is a Galois-equivariant isogeny
\[ \varphi_A:T_A(\bar\Q_p)/L_A\longrightarrow A(\bar\Q_p) \] 
defined over $K_p$ (Theorem \ref{homothetic-thm}). At least when the reduction of $A$ at $p$ is split (a condition that we do not assume in the main body of the article), an explicit description of $\varphi_A$ can be given in terms of the Galois-equivariant analytic isomorphism
\[ (\bar\Q_p^\times)^d\big/\big\langle\q_1,\dots,\q_d\big\rangle\overset\simeq\longrightarrow A(\bar\Q_p) \]
over $\Q_p$ due (among others) to Tate, Morikawa and Mumford. Here $\langle\q_1,\dots,\q_d\rangle$ is the lattice inside $\Q_p^d$ generated by the Tate periods $\q_1,\dots,\q_d$ for $A$ at $p$, which turns out to be homothetic to $L_A$ (Theorem \ref{homothetic2-thm}). The reader can find details about these constructions in Appendix \ref{appendix}.
 
Regarding $\epsilon$ as fixed, for every optimal embedding $\psi$ of $\cO_c$ into $R_0$ we define
\[ P_{A,\psi}:=\varphi_A\bigl(\mathcal P_{\mathbb H^\epsilon_A,\psi}\bigr)\in A(K_p). \]
These are the Darmon points (of conductor $c$) on $A$ alluded to before. To state our formula for the points $P_{A,\psi}$, let $\C_p$ be the completion of $\bar\Q_p$ and recall that the theory of Lie groups gives a logarithm map
\[ \log_A:A(\C_p)\longrightarrow\text{Lie}(A(\C_p))\simeq\C_p^d \]
on the $p$-adic points of $A$. In particular, $\log_A(P)\in K_p^d$ for all $P\in A(K_p)$. Define $\Psi_A:=\log_A\circ\,\varphi_A$. Now let $f$ be an endomorphism of $A_{/\C_p}$ and let $f_*:\text{Lie}(A(\C_p))\rightarrow\text{Lie}(A(\C_p))$ be the linear map induced by $f$. 

The second main result of this paper, which corresponds to Theorem \ref{thm2} and can be deduced from Theorem \ref{intro-thm1} by taking $\mathbb H=\mathbb H_A^\epsilon$ and $\Psi=\Psi_A$, is 

\begin{theorem} \label{intro-thm2}
$\displaystyle{\log_A\bigl(f(P_{A,\psi})\bigr)=-t\cdot f_*\left(\int_\X\Psi_A(x-z_\psi y)d{\tilde\mu}_{A,\gamma_\psi}(x,y)\right)}$. 
\end{theorem}

By choosing $f=\text{id}_A$ in Theorem \ref{intro-thm2}, we obtain

\begin{corollary} \label{intro-coro}
$\displaystyle{\log_A(P_{A,\psi})=-t\cdot\int_\X\Psi_A(x-z_\psi y)d{\tilde\mu}_{A,\gamma_\psi}(x,y)}$.
\end{corollary}

This is the first formula for points of Stark--Heegner type on higher-dimensional abelian varieties. When $A$ is an elliptic curve Corollary \ref{intro-coro} essentially reduces to \cite[Corollary 3.4]{LV-darmon}. However, it might be useful to notice that 
\begin{itemize}
\item the result in \cite[Corollary 3.4]{LV-darmon} only deals with Darmon points of conductor $c=1$, while in Theorem \ref{intro-thm2} there are no restrictions on $c$ (except for the standard ones requiring $c$ to be coprime with $MDd_Kp$); 
\item Corollary \ref{intro-coro} expresses the logarithm of $P_{A,\psi}$ in terms of integrals against an \emph{arbitrary} lift $\boldsymbol{\tilde\mu}_A$, while in \cite[Corollary 3.4]{LV-darmon} a specific choice is made (using a control theorem for cohomology groups in Hida families). 
\end{itemize} 
The above restrictions in the case of elliptic curves are crucial for the proof of the rationality results in \cite{LV-darmon}, but the more general approach presented here does not require any of these assumptions. 

The interest in having such formulas at our disposal is (at least) twofold. On the one hand, since the construction of measure-valued cocycles in \cite{LRV1} is entirely explicit, Theorem \ref{intro-thm2} allows us to compute (at least in principle) the points in the orbit of any Darmon point on $A$ under the action of the endomorphism ring of $A$ (cf. also Corollary \ref{thm2-coro}). On the other hand, an integral expression for the logarithm of $P_{A,\psi}$ like the one of Corollary \ref{intro-coro} will presumably play a crucial role in the proof of rationality results for genus character combinations of Darmon points on $A$ of the type obtained in \cite{LV-darmon} in the case of elliptic curves. We plan to investigate these questions in a future project.

\section{Measure-valued cohomology} \label{sec-2.1} \label{sec-2.2} 

\subsection{Homology of Shimura curves and measures}

We recall notation and results of \cite{LRV1}. Let $D>1$ be a square-free product of an \emph{even} number of primes, let $M\geq1$ be an integer prime to $D$ and fix a prime number $p$ not dividing $MD$. Let $B$ be the (indefinite) quaternion algebra over $\Q$ of discriminant $D$. Let $R_0\subset R_1$ be Eichler orders of $B$ of level $M$ and $Mp$, respectively. For $i=0,1$ write $\Gamma_i$ for the group of units of norm $1$ in $R_i$ and consider the (compact) Shimura curve $X_i:=\Gamma_i\backslash\mathcal H$, where $\mathcal H$ is the complex upper half-plane and the action of $B^\times$ on $\PP^1(\C)$ is by M\"obius (i.e., fractional linear) transformations via a fixed isomorphism of $\R$-algebras
\[ i_\infty:B\otimes_\Q\R\overset\simeq\longrightarrow\M_2(\R). \]
Let $\om_p$ be an element in $R_0$ of reduced norm $p$ that normalizes $\Gamma_1$ and denote 
\[ \pi_1,\pi_2: X_1\longrightarrow X_0,\qquad\Gamma_1z\overset{\pi_1}\longmapsto\Gamma_0 z,\qquad\Gamma_1z\overset{\pi_2}\longmapsto\Gamma_0\om_p z \]
the two degeneracy maps. Moreover, let
\[ \pi^\ast:=\pi_1^\ast\oplus\pi_2^\ast:H_1(X_0,\Z)^2\longrightarrow H_1(X_1,\Z) \]
be the map induced in singular homology by pull-back. In terms of group homology, it corresponds to the map
\[ H_1(\Gamma_0,\Z)^2\longrightarrow H_1(\Gamma_1,\Z)\simeq\Gamma_1^{\rm ab} \]
induced by corestriction. Let $H$ denote the maximal torsion-free quotient of the cokernel of $\pi^\ast$. If $J_1^{\text{$p$-{\rm new}}}$ is the $p$-new quotient (i.e., the maximal toric quotient at $p$) of the Jacobian variety $J_1$ of $X_1$ and $g$ is the dimension of $J_1^{\text{$p$-{\rm new}}}$ then $H$ is a free abelian group of rank $2g$. Throughout this paper fix a non-zero torsion-free quotient $\mathbb H$ of $H$. 

Fix an isomorphism of $\Q_p$-algebras
\[ i_p:B\otimes_\Q\Q_p\overset\simeq\longrightarrow\M_2(\Q_p) \] 
such that $i_p(R_0\otimes\Z_p)=\M_2(\Z_p)$ and $i_p(R_1\otimes\Z_p)$ is the order of $\M_2(\Q_p)$ consisting of the matrices $\smallmat abcd\in\M_2(\Z_p)$ with $c\equiv0\pmod{p}$. Denoting $\C_p$ the completion of an algebraic closure of $\Q_p$, let $B^\times$ act on $\PP^1(\C_p)$ by fractional linear transformations via $i_p$. 

Let ${\rm n}:B\rightarrow\Q$ denote the reduced norm map and define Ihara's group $\Gamma$ as 
\begin{equation} \label{Gamma}
\Gamma:=\bigl\{\gamma\in R_0\otimes\Z[1/p]\mid{\rm n}(\gamma)=1\bigr\}\;\overset{i_p}\longmono\;\SL_2(\Q_p).
\end{equation}
Write $\mathcal M$ for the group of $\mathbb H$-valued measures on $\PP^1(\Q_p)$ and $\mathcal M_0$ for the subgroup of $\mathcal M$ consisting of those measures  with total mass $0$. There is a canonical left action of $\GL_2(\Q_p)$ on $\mathcal M$ and $\mathcal M_0$ given by the integration formula 
\begin{equation} \label{int-formula}
\int_{\PP^1(\Q_p)}\varphi(t)d(\gamma\cdot\nu)(t):=\int_{\PP^1(\Q_p)}\varphi\left(\frac{at+b}{ct+d}\right)d\nu 
\end{equation}
for all step functions $\varphi:\PP^1(\Q_p)\rightarrow\C_p$ and all $\gamma=\smallmat abcd\in\GL_2(\Q_p)$. Then $B^\times$ acts on $\mathcal M$ and $\mathcal M_0$ via $i_p$ and the embedding $B\hookrightarrow B\otimes\Q_p$. 

\subsection{Harmonic cocycles, radial systems and Hecke algebras}

In the following lines we review the construction of a canonical element 
\[ \boldsymbol\mu_{\mathbb H}\in H^1(\Gamma,\mathcal M_0) \]
introduced in \cite[Section 4]{LRV1}. 

Let $\mathcal T$ be the Bruhat--Tits tree of $\PGL_2(\Q_p)$ (see \cite[Ch. II, \S 1]{Se}), whose set of vertices (respectively, oriented edges) will be denoted $\mathcal V$ (respectively, $\mathcal E$). The unit groups $B^\times$ and $(B\otimes\Q_p)^\times$ act on the left on $\mathcal T$ via $i_p$. Let $v_*\in\mathcal V$ denote the distinguished vertex of $\mathcal T$ corresponding to the maximal order $\M_2(\Z_p)$ of $\M_2(\Q_p)$, and say that a vertex $v\in\mathcal V$ is \emph{even} (respectively, \emph{odd}) if its distance from $v_*$ is even (respectively, odd). Moreover, denote $\hat v_*$ the vertex corresponding to $\smallmat{\Z_p}{p^{-1}\Z_p}{p\Z_p}{\Z_p}$, fix the edge $e_*=(v_*,\hat v_*)\in\mathcal E$ and say that an edge $e=(v_1,v_2)\in\mathcal E$ is \emph{even} (respectively, \emph{odd}) if $v_1$ is even (respectively, odd). Write $\mathcal E^+$ for the set of even vertices of $\mathcal T$. Finally, if $e=(v_1,v_2)$ write $\bar e$ for the reversed edge $(v_2,v_1)$.

For any vertex $v\in\mathcal V$ define 
\[ s(v):=\sum_{e=(v,v')}c(e). \] 
If $A$ is an abelian group then a function $c:\mathcal E\rightarrow A$ is said to be an \emph{harmonic cocycle} with values in $A$ if
\begin{itemize} 
\item $c(\bar e)=-c(e)$ for all $e\in\mathcal E$,
\item $s(v)=0$ for all $v\in\mathcal V$.
\end{itemize} 
Let $\mathcal F_{\rm har}$ denote the abelian group of $\mathbb H$-valued harmonic cocycles. There is a canonical isomorphism
\begin{equation} \label{HC-Meas}
H^1(\Gamma,\mathcal F_{\rm har})\overset\simeq\longrightarrow H^1(\Gamma,\mathcal M_0)
\end{equation}
given by the following recipe. Recall that $\Gamma$ acts transitively on $\mathcal E$ and the stabilizer of $e_*$ in $\Gamma$ is $\Gamma_1$. For any $e\in\mathcal E$ choose $\gamma\in\Gamma$ such that $\gamma(e_*)=e$ and set $U_e:=\gamma(\Z_p)$. The family $\{U_e\}_{e\in\mathcal E}$ is a basis of compact open subsets of $\PP^1(\Q_p)$. 
Given any $c\in\mathcal F_{\rm har}$, we may define a measure $\nu\in\mathcal M_0$ by the formula $\nu(U_e):=c(e)$. This sets up a canonical isomorphism between 
$\mathcal F_{\rm har}$ and $\mathcal M_0$, from which one deduces \eqref{HC-Meas}. For details see, e.g., \cite[Lemma 27]{Gr} or \cite[\S 4.1]{LRV1}. 

In light of \eqref{HC-Meas}, our task is to describe a canonical element in $H^1(\Gamma,\mathcal F_{\rm har})$. Set $\hat\Gamma_0:=\omega_p\Gamma_0\omega_p^{-1}$. A  system of representatives $\mathcal Y=\{\gamma_e\}_{e\in\mathcal E^+}$ for the cosets in $\Gamma_1\backslash\Gamma$ is said to be \emph{radial} if there are sets of representatives
\begin{itemize}
\item $\{\gamma_i\}_{i=0}^p$ for the cosets in $\Gamma_1\backslash \Gamma_0$,
\item $\{\tilde\gamma_i\}_{i=0}^p$ for the cosets in $\Gamma_1\backslash\hat\Gamma_0$,
\item $\{\gamma_v\}_{v\in\mathcal V^+}$ for the cosets in $\Gamma_0\backslash\Gamma$,
\item $\{\gamma_v\}_{v\in\mathcal V^-}$ for the cosets in $\hat\Gamma_0\backslash\Gamma$  
\end{itemize}
such that the following conditions hold: 
\begin{itemize}
\item $\{\gamma_e\}_{e=(v,v')}=\{\gamma_i\gamma_v\}_{i=0}^{p}$ for all $v\in\mathcal V^+$, 
\item $\{\gamma_e\}_{e=(v',v)}=\{\tilde\gamma_i\gamma_v\}_{i=0}^{p}$ for all $v\in\mathcal V^-$,
\item $\gamma_0=\tilde\gamma_0=\gamma_{v_*}=\gamma_{\hat v_*}=1$. 
\end{itemize}  
Choose a radial system of representatives $\mathcal Y$ for the cosets in ${\Gamma_1}\backslash\Gamma$ (which exists by \cite[Proposition 4.8]{LRV1}) and define $\mathcal F_0$ to be the abelian group of $\mathbb H$-valued functions $\varphi$ on $\mathcal E$ such that $\varphi(\bar e)=-\varphi(e)$. The group $\mathcal F_0$ is endowed with the action of $B^\times$ given by the formula $(g\cdot \varphi)(e):=\varphi(g^{-1}(e))$, which agrees with the action already defined on the subgroup $\mathcal F_{\rm har}$ of $\mathcal F_0$ via \eqref{int-formula} and isomorphism \eqref{HC-Meas}. 

Keeping in mind that $\mathbb H$ is (canonically isomorphic to) a quotient of ${\Gamma_1}$, define, as in \cite[Definition 4.2]{LRV1}, the function 
\[ \mu_{\mathbb H}^\mathcal Y:\Gamma\longrightarrow\mathcal F_0 \] 
associated with $\mathcal Y$ and $\mathbb H$ by the following rules:
\begin{itemize} 
\item if $\gamma\in\Gamma$ and $e\in\mathcal E^+$ let $g_{\gamma,e}\in{\Gamma_1}$ be defined by the equation $\gamma_e\gamma=g_{\gamma,e}\gamma_{\gamma^{-1}(e)}$, then set $\mu_{\mathbb H}^\mathcal Y(\gamma)(e):=[g_{\gamma,e}]\in\mathbb H$; 
\item if $\gamma\in\Gamma$ and $e\not\in\mathcal E^+$ then set $\mu_{\mathbb H}^\mathcal Y(e):=-\mu_{\mathbb H}^\mathcal Y(\bar e)$. 
\end{itemize} 
One can show that $\mu^\mathcal Y_{\mathbb H}$ is a $1$-cocycle on $\Gamma$ with values in $\mathcal F_{\rm har}$ (see \cite[Proposition 4.8]{LRV1}). Furthermore, let \[ \varrho:H^1(\Gamma,\mathcal F_{\rm har})\longrightarrow H^1(\Gamma,\mathcal F_0) \]  
denote the map arising by taking $\Gamma$-cohomology of the short exact sequence
\[ 0\longrightarrow\mathcal F_{\rm har}\longrightarrow\mathcal F_0\longrightarrow\mathcal F_1\longrightarrow 0 \]
where $\mathcal F_1$ denotes the abelian group of $\mathbb H$-valued functions on $\mathcal V$. Then it can be easily checked that the image via $\varrho$ of the class represented by $\mu_{\mathbb H}^\mathcal Y$ does not depend on the choice of $\mathcal Y$ (\cite[Proposition 4.3]{LRV1}). However, we need to slightly modify $\mu_{\mathbb H}^\mathcal Y$ in order to make it independent of the choice of $\mathcal Y$ as an element of $H^1(\Gamma,\mathcal F_{\rm har})$ (see the discussions in \cite[Remarks 4.5 and 4.9]{LRV1}). To do this, we briefly recall the Hecke action on this cohomology group.  

For a prime $\ell\nmid Mp$ let $S_\ell$ denote the set of elements in $R_0\otimes_\Z\Z_\ell$ with non-zero norm. For a prime $\ell|Mp$ let $n_\ell$ be the maximal power of $\ell$ dividing $Mp$. For primes $\ell\neq p$ with $\ell|M$ fix an isomorphism of $\Q_\ell$-algebras 
\[ i_\ell:B\otimes_\Q\Q_\ell\overset\simeq\longrightarrow\M_2(\Q_\ell) \] 
such that $i_\ell(R_0\otimes\Z_\ell)$ is the order consisting of matrices $\smallmat abcd\in\M_2(\Z_\ell)$ with $\ell^{n_\ell}|c$. For all primes $\ell|Mp$ define $S_\ell$ to be the inverse image via $i_\ell$ of the semigroup consisting of matrices $g=\smallmat abcd\in\M_2(\Z_\ell)$ such that $\ell^{n_\ell}|c$, $a\in\Z_\ell^\times$ and $\det(g)\neq0$. Denote $S$ the subsemigroup of $B^\times$ defined as in \cite[\S 2.3]{LRV1} (where it is denoted $S_{(p,M)}$) by 
\begin{equation} \label{S}
S:=B^\times\cap\prod_{\ell\neq p}S_\ell,
\end{equation} 
the product being taken over all primes $\ell\not=p$. Finally, let $\T$ stand for the Hecke algebra associated with the pair $(\Gamma,S)$ (see \cite[\S 2.3]{LRV1} for details). As recalled in \cite[\S 2.1]{LRV1}, the algebra $\T$ acts on the cohomology group $H^1(\Gamma,M)$ for any left $S$-module $M$. In particular, for every prime $r$ we have an action of the Hecke operator $T_r\in\T$ on $H^1(\Gamma,\mathcal F_0)$. 

Now define 
\[ \boldsymbol\mu_{\mathbb H}\in H^1(\Gamma,\mathcal F_{\rm har}) \] 
to be the class obtained by fixing a prime number $r\nmid MDp$ and applying the Hecke operator $t_r:=T_r-(r+1)\in\T$ to the class of $\mu_{\mathbb H}^\mathcal Y$ in $H^1(\Gamma,\mathcal F_{\rm har})$. As shown in \cite[Lemma 4.11]{LRV1}, this operation has the effect of making $\boldsymbol\mu_{\mathbb H}$ independent of the choice of $\mathcal Y$ made above. Finally, fix a representative $\mu_{\mathbb H}$ of $\boldsymbol\mu_{\mathbb H}$. We will use the same symbol $\boldsymbol\mu_{\mathbb H}$ for the element in $H^1(\Gamma,\mathcal M_0)$ corresponding to $\boldsymbol\mu_{\mathbb H}$ via isomorphism \eqref{HC-Meas}. 
 
\section{Darmon points on $p$-adic tori} \label{sec2.5}

\subsection{Homology and $p$-adic integration}

Set 
\[ T_{\mathbb H}:=\mathbb G_m\otimes_\Z\mathbb H \]
where $\mathbb G_m$ is the multiplicative group (viewed as a functor on commutative $\Q$-algebras).  

Write $\mathcal H_p:=\C_p-\Q_p$ for Drinfeld's $p$-adic plane, let $D:=\Div(\mathcal H_p)$ be the group of divisors on $\mathcal H_p$ and let $D_0:=\Div^0(\mathcal H_p)$ be the subgroup of divisors of degree $0$. Then there is a short exact sequence 
\begin{equation} \label{ShExSq}
0\longrightarrow D_0\longrightarrow D\xrightarrow{\rm deg}\Z\longrightarrow0 
\end{equation}
induced by the degree map. The long exact sequence in homology associated with \eqref{ShExSq} gives a boundary map 
\begin{equation} \label{partial}
\partial:H_2(\Gamma,\Z)\longrightarrow H_1(\Gamma,D_0).
\end{equation}
For any $d\in D_0$ choose a rational function $f_d$ on $\PP^1(\C_p)$ with divisor $d$, and for all $\nu\in\mathcal M_0$ define 
\[ \langle d,\nu\rangle=\mint_{\PP^1(\Q_p)}f_dd\nu:=\lim_{\|\mathcal U\|\rightarrow0}\prod_{U\in\mathcal U}f_d(t_U)\otimes\nu(U)\in T_{\mathbb H}(\C_p). \] 
Here the limit is taken over increasingly fine covers $\mathcal U$ of $\PP^1(\Q_p)$ by disjoint compact open subsets and $t_U\in U$ is a sample point. This limit of Riemann products converges because $\nu$ is a measure, and is independent of the choice of $f_d$ because $\nu$ has total mass $0$. The above formula sets up a pairing 
\[ \langle\,,\rangle:D_0\times\mathcal M_0\longrightarrow T_{\mathbb H}(\C_p). \] 
By construction, this pairing factors through $H_0(\Gamma,D_0\otimes_\Z\mathcal M_0)$ and thus, by cap product, we can also define a pairing 
\begin{equation} \label{pairing}
H_1(\Gamma,D_0)\times H^1(\Gamma,\mathcal M_0)\longrightarrow T_{\mathbb H}(\C_p).
\end{equation} 
By fixing $\boldsymbol\mu_{\mathbb H}\in H^1(\Gamma,\mathcal M_0)$ in \eqref{pairing} we get a map 
\begin{equation} \label{I}
\mathcal I:H_1(\Gamma,D_0)\longrightarrow T_{\mathbb H}(\C_p)
\end{equation}
canonically attached to $\boldsymbol\mu_{\mathbb H}$. Finally, composing \eqref{partial} and \eqref{I} we obtain a map
\[ \Phi:=\mathcal I\circ\partial:H_2(\Gamma,\Z)\longrightarrow T_{\mathbb H}(\C_p). \]
Denote $L_{\mathbb H}$ the image of $\Phi$. One can show that $L_{\mathbb H}$ is a Hecke-stable lattice contained in $T_{\mathbb H}(\Q_p)$ (see \cite[Proposition 6.1]{LRV1}). 

Fix $z\in K_p-\Q_p$ and let $\boldsymbol{d}_z\in H^2\bigl(\Gamma,T_{\mathbb H}(\mathbb C_p)\bigr)$ be the cohomology class represented by the $2$-cocycle 
\begin{equation} \label{d-cocycle-eq}
d_z:\Gamma\times\Gamma\longrightarrow T_{\mathbb H}(K_p),\qquad(\gamma_1,\gamma_2)\longmapsto\mint_{\PP^1(\Q_p)}\frac{s-\gamma_1^{-1}(z)}{s-z}d\mu_{\mathbb H,\gamma_2}(s). 
\end{equation}
The class $\boldsymbol{d}_z$ does not depend on the choice of the representative $\mu_{\mathbb H}$ of $\boldsymbol\mu_{\mathbb H}$. Write $\bar d_z$ for the composition of $d_z$ with the canonical projection onto $T_{\mathbb H}(K_p)/L_{\mathbb H}$ and denote $\boldsymbol{\bar{d}}_z$ the resulting class. By construction, $L_{\mathbb H}$ is the smallest subgroup of $T_{\mathbb H}(\Q_p)$ such that $\boldsymbol{\bar{d}}_z$ is trivial in $H^2(\Gamma,T_{\mathbb H}/L_{\mathbb H})$, so there exists 
\[ \beta_z:\Gamma\longrightarrow T_{\mathbb H}/L_{\mathbb H} \] 
such that 
\[ \bar d_z(\gamma_1,\gamma_2)=\beta_z(\gamma_1\gamma_2)\cdot\beta_z(\gamma_1)^{-1}\cdot\beta_z(\gamma_2)^{-1} \]
for all $\g_1,\g_2\in\Gamma$. The map $\beta_z$ is well defined only up to elements in $\Hom(\Gamma,T_{\mathbb H}/L_{\mathbb H})$. To deal with this ambiguity, recall that the abelianization $\Gamma^{\rm ab}$ of $\Gamma$ is finite (\cite[Proposition 2.1]{LRV2}). Hence if $t$ is the exponent of $\Gamma^{\rm ab}$ then $t\cdot\beta_z$ is well defined. However, since $d_z$ depends on the choice of $\mu_{\mathbb H}$, the map $t\cdot\beta_z$ may depend on the choice of a representative for $\boldsymbol\mu_{\mathbb H}$.

\subsection{Darmon points on $p$-adic tori} 

Let $K=\Q(\sqrt{d_K})$ be a real quadratic field with discriminant $d_K$ such that all primes dividing $Dp$ (respectively, $M$) are inert (respectively, split) in $K$, and fix an embedding $K\hookrightarrow\R$. Let $c\geq1$ be an integer prime to $MDd_Kp$, let $\cO_c:=\Z+c\cO_K$ be the order of $K$ of conductor $c$, let $H_c^+$ be the narrow ring class field of $K$ of conductor $c$ and set $G^+_c:=\Gal(H_c^+/K)$. The group $G_c^+$ is isomorphic to the narrow class group of $\cO_c$ via the reciprocity map of global class field theory.

Let $\psi:K\hookrightarrow B$ be an optimal embedding of $\cO_c$ into $R_0$, i.e. an embedding of $K$ into $B$ such that 
\[\psi(\cO_c)=\psi(K)\cap R_0.\] 
Write $\Emb(\cO_c,R_0)$ for the set of such embeddings. The group $\Gamma_0$ acts on $\Emb(\cO_c,R_0)$ by conjugation.

Let $\tau$ denote the generator of $\Gal(K_p/\Q_p)$. If $\psi:K\hookrightarrow B$ is an embedding of $\Q$-algebras then the quadratic form $Q_\psi(x,y)$ with coefficients in $\Q_p$ associated with $\psi$ is defined by 
\[ Q_\psi(x,y):=cx^2-2axy-by^2\qquad\text{where $i_p\bigl(\psi(\sqrt{d_K})\bigr)=\begin{pmatrix}a&b\\c&{-a}\end{pmatrix}$}. \] 
We can factor $Q_\psi(x,y)$ as 
\[ Q_\psi(x,y)=c(x-z_\psi y)(x-\bar z_\psi y) \] 
where $z_\psi,\bar z_\psi\in K_p-\Q_p$ are the roots of the equation $cz^2-2az-b=0$ and $\tau(z_\psi)=\bar z_\psi$. The two roots $z_\psi$ and $\bar z_\psi$ are the only fixed points for the action of $\psi(K^\times)$ on $\PP^1(K_p)$ by fractional linear transformations via $i_p$. We may also order $z_\psi$ and $\bar z_\psi$ by requiring that  
\[ \psi(\alpha)\binom{z_\psi}{1}=\alpha\binom{z_\psi}{1},\qquad\psi(\alpha)\binom{\bar z_\psi}{1}=\tau(\alpha)\binom{\bar z_\psi}{1} \]
for all $\alpha\in K$.

By Dirichlet's unit theorem, the abelian group of units of norm $1$ in $\cO_c$ is free of rank $1$. Let $\varepsilon_c$ be a generator of this group such that $\varepsilon_c>1$ with respect to the fixed embedding $K\hookrightarrow\R$ and set $\gamma_\psi:=\psi(\varepsilon_c)\in{\Gamma_0}$.

\begin{definition} \label{darmon-def}
The \emph{Darmon points of conductor $c$} on $T_{\mathbb H}(K_p)/L_{\mathbb H}$ are the points 
\[ \mathcal P_{\mathbb H,\psi}:=t\cdot\beta_{z_\psi}(\gamma_\psi)\in T_{\mathbb H}(K_p)/L_{\mathbb H} \]
where $\psi$ varies in $\Emb(\cO_c,R_0)$.
\end{definition}

The points $\mathcal P_{\mathbb H,\psi}$ (or, rather, their images on abelian varieties) are expected to be rational over $H_c^+$ and to satisfy a suitable Shimura reciprocity law under the action of $G^+_c$: the reader is referred to \cite[\S 3.2]{LRV2} for precise conjectures and to \cite{LV-darmon} for partial results in this direction.

The next two propositions generalize to $\mathbb H$ results that were proved in \cite{LRV2} for $H$. Although the proofs are similar, we add some details for the convenience of the reader. 

\begin{proposition} \label{lemma2.2}
The point $\mathcal P_{\mathbb H,\psi}$ does not depend on the choice of a representative of $\boldsymbol\mu_{\mathbb H}$. 
\end{proposition}

\begin{proof} Let $\mu$ and $\mu'$ be two representatives for $\boldsymbol\mu_{\mathbb H}$, so that there exists $m\in\mathcal M_0$ such that $\mu_\g=\mu'_\g+\g(m)-m$ for all $\g\in\Gamma$. Then  
\[ \begin{split}
   d_{z_\psi}(\g_1,\g_2)\div d'_{z_\psi}(\g_1,\g_2)=&\mint_{\PP^1(\Q_p)}\frac{t-\gamma_1^{-1}z}{t-z}d(\gamma_2m-m)\\
   =&\mint_{\PP^1(\Q_p)}\frac{t-\gamma_2^{-1}\gamma_1^{-1}z}{t-\gamma_2^{-1}z}dm\div\mint_{\PP^1(\Q_p)}\frac{t-\gamma_1^{-1}z}{t-z}dm\\
   =&\mint_{\PP^1(\Q_p)}\frac{t-\gamma_2^{-1}\gamma_1^{-1}z}{t-z}dm\\
   &\div\bigg(\mint_{\PP^1(\Q_p)}\frac{t-\gamma_1^{-1}z}{t-z}dm\cdot\mint_{\PP^1(\Q_p)}\frac{t-\gamma_2^{-1}z}{t-z}dm\bigg).
   \end{split} \] 
If we set 
\begin{equation} \label{eq-k}
\nu(\g):=\mint_{\PP^1(\Q_p)}\frac{s-\gamma^{-1}(z_\psi)}{s-z_\psi}dm(s)
\end{equation}
then 
\[ d_{z_\psi}(\g_1,\g_2)=d'_{z_\psi}(\g_1,\g_2)\cdot\nu(\g_1\g_2)\cdot\nu(\g_1)^{-1}\cdot\nu(\g_2)^{-1} \]
for all $\g_1,\g_2\in\Gamma$. (In particular, this shows that the $2$-cocycles $d_{z_\psi}$ and $d'_{z_\psi}$ defined as in \eqref{d-cocycle-eq} by using $\mu$ and $\mu'$, respectively, are cohomologous.) Write $\bar\nu$ for the composition of $\nu$ with the projection onto $T_{\mathbb H}(K_p)/L_{\mathbb H}$. If $\beta_{z_\psi}:\Gamma\rightarrow T_{\mathbb H}/L_{\mathbb H}$ (respectively, $\beta'_{z_\psi}:\Gamma\rightarrow T_{\mathbb H}/L_{\mathbb H}$) splits $d_{z_\psi}$ (respectively, $d'_{z_\psi}$) modulo $L_{\mathbb H}$ then 
\[ \beta_{z_\psi}=\beta'_{z_\psi}\cdot\bar\nu\cdot\varphi' \] 
for a suitable homomorphism $\varphi':\Gamma\rightarrow T_{\mathbb H}/L_{\mathbb H}$. It follows that 
\begin{equation} \label{eq--}
t\cdot\beta_{z_\psi}=\bigl(t\cdot\beta'_{z_\psi}\bigr)\cdot(t\cdot\bar\nu).
\end{equation} 
Since $\gamma_\psi(z_\psi)=z_\psi$ and $m$ has total mass $0$, equation \eqref{eq-k} shows that $t\cdot\nu(\gamma_\psi)=1$. By definition of the Darmon point $\mathcal P_{\mathbb H,\psi}$, the result follows from this and equation \eqref{eq--}. \end{proof}

The following proposition studies the dependence of $\mathcal P_{\mathbb H,\psi}$ on $\psi$.
 
\begin{proposition} \label{conjugacy-prop}
The point $\mathcal P_{\mathbb H,\psi}$ depends only on the $\Gamma_0$-conjugacy class of $\psi$. 
\end{proposition}

\begin{proof} Fix $\gamma\in\Gamma_0$ and set $\psi':=\gamma\psi\gamma^{-1}$. Recall the radial system $\mathcal Y=\{\gamma_e\}_{e\in\mathcal E^+}$ used to compute $\boldsymbol\mu_{\mathbb H}$ and introduce the set 
\[ \mathcal Y':=\bigl\{\gamma'_e:=\gamma\gamma_{\gamma^{-1}(e)}\gamma^{-1}\bigr\}_{e\in\mathcal E^+}. \] 
One can check that $\mathcal Y'$ is again a radial system. For this, define $\gamma_i':=\gamma\gamma_i\gamma^{-1}$, $\tilde\gamma_i':=\gamma\tilde\gamma_i\gamma^{-1}$ and $\gamma_{\gamma(v)}':=\gamma\gamma_v\gamma^{-1}$. Also, if $e=(v,v')$ then write $s(e):=v$ and $t(e):=v'$. It follows that 
\[ \bigl\{\gamma_i'\gamma_{\gamma(v)}'\bigr\}_{i=0}^p=\bigl\{\gamma\gamma_i\gamma_{\gamma(v)}\gamma^{-1}\bigr\}_{i=0}^p=\bigl\{\gamma\gamma_e\gamma^{-1}\bigr\}_{s(e)=\gamma(v)}=\bigl\{\gamma_{\gamma(e)}'\bigr\}_{s(\gamma(e))=\gamma(v)} \]
for $v\in\mathcal V^+$, and
\[ \bigl\{\tilde\gamma_i'\gamma_{\gamma(v)}'\bigr\}_{i=0}^p=\bigl\{\gamma\tilde\gamma_i\gamma_{\gamma(v)}\gamma^{-1}\bigr\}_{i=0}^p=\bigl\{\gamma\gamma_e\gamma^{-1}\bigr\}_{s(e)=\gamma(v)}=\bigl\{\gamma_{\gamma(e)}'\bigr\}_{s(\gamma(e))=\gamma(v)} \]
for $v\in\mathcal V^-$. To prove the last equalities in the above equations, observe that  
\[ \bigl\{e\mid s(e)=\gamma(v)\bigr\}=\bigl\{\gamma(e)\mid s(\gamma(e))=\gamma(v)\bigr\} \] 
and that the same is true for $t(e)$ replacing $s(e)$. Since $\gamma\in\Gamma_0$, one has that
\[ \gamma(v)\in\begin{cases}\mathcal V^+& \text{if $v\in\mathcal V^+$},\\\mathcal V^-& \text{if $v\in\mathcal V^-$}. \end{cases} \]
But $\Gamma_0$ acts transitively on the sets $\mathcal V^+$ and $\mathcal V^-$, hence 
\[ \bigl\{\gamma'_e\bigr\}_{s(e)=v}=\bigl\{\gamma_i'\gamma_{\gamma(v)}'\bigr\}_{i=0}^p\quad\text{for $v\in\mathcal V^+$},\qquad\bigl\{\gamma'_e\bigr\}_{t(e)=v}=\bigl\{\tilde\gamma_i'\gamma_{\gamma(v)}'\bigr\}_{i=0}^p\quad\text{for $v\in\mathcal V^-$}. \]
The final observation to prove that $\mathcal Y'$ is radial is that, trivially, $\gamma'_0=\tilde\gamma'_0=\gamma'_{v_*}=\gamma'_{\hat v_*}=1$. 

By definition, $\mu_{\rm univ}^{\mathcal Y_{\rm rad}}(\gamma_2)(e)=[g_{\gamma_2,e}]$ with $\gamma_e\gamma_2=g_{\gamma_2,e}\gamma_{e'}$. Conjugating by $\gamma$ gives
\[ \gamma\gamma_e\gamma^{-1}\gamma\gamma_2\gamma^{-1}=\gamma g_{\gamma_2,e}\gamma^{-1}\gamma^{-1}\gamma\gamma_{e'}\gamma^{-1}, \]
whence $\gamma'_{\gamma(e)}\gamma\gamma_2\gamma^{-1}=\gamma g_{\gamma_2,e}\gamma^{-1}\gamma'_{\gamma(e')}$. This shows that 
\[ \mu^{\mathcal Y'}_{\mathbb H}(\gamma\gamma_2\gamma^{-1})\bigl(\gamma(e)\bigr)=\bigl[\gamma g_{\gamma_2,e}\gamma^{-1}\bigr]=[g_{\gamma_2,e}]=\mu_{\mathbb H}^{\mathcal Y}(\gamma_2)(e). \]  
Applying $t_r$ to the above equation, we get 
\begin{equation} \label{eq1}
\gamma^{-1}\mu^{\mathcal Y'}_{\mathbb H}(\gamma\gamma_2\gamma^{-1})(e)=\mu^{\mathcal Y}_{\mathbb H}(\gamma_2)(e).
\end{equation} 
Therefore
\begin{equation} \label{eq-prop2.4}
\begin{split}
\mint_{\PP^1(\Q_p)}\frac{s-\gamma_1^{-1}z_\psi}{s-z_\psi}d\mu^{\mathcal Y}_{\mathbb H,\gamma_2}(s)
&=\mint_{\PP^1(\Q_p)}\frac{s-\gamma_1^{-1}z_\psi}{s-z_\psi}\gamma^{-1}d\mu^{\mathcal Y'}_{\mathbb H,\gamma\gamma_2\gamma^{-1}}(s)\\
&=\mint_{\PP^1(\Q_p)}\frac{s-\gamma\gamma_1^{-1}z_\psi}{s-\gamma z_\psi}d\mu^{\mathcal Y'}_{\mathbb H,\gamma\gamma_2\gamma^{-1}}(s)\\
&=\mint_{\PP^1(\Q_p)}\frac{s-\gamma\gamma_1^{-1}\gamma^{-1}\gamma z_\psi}{s-\gamma z_\psi}d\mu^{\mathcal Y'}_{\mathbb H,\gamma\gamma_2\gamma^{-1}}(s)\\
&=\mint_{\PP^1(\Q_p)}\frac{s-\gamma\gamma_1^{-1}\gamma^{-1}z_{\gamma\psi\gamma^{-1}}}{s-z_{\gamma\psi\gamma^{-1}}}d\mu^{\mathcal Y'}_{\mathbb H,\gamma\gamma_2\gamma^{-1}}(s).
\end{split}
\end{equation} 
The first equality follows from \eqref{eq1}, the second follows because $\langle\gamma\cdot f,\gamma\cdot\nu\rangle=\langle f,\nu\rangle$ and, if $d$ is a divisor, $\gamma\cdot f_d=f_{\gamma d}$ (for these properties, see \cite[\S 5.1]{LRV1}), the third is obvious and the fourth is due to the fact that $\gamma z_\psi=z_{\gamma\psi\gamma^{-1}}$. If $\beta'_{z_{\psi'}}$ splits the $2$-cocycle 
\[ (\gamma_1,\gamma_2)\longmapsto\mint_{\PP^1(\Q_p)}\frac{s-\gamma_1^{-1}z_{\psi'}}{s-z_{\psi'}}d\mu^{\mathcal Y_{\rm rad}'}_{\mathbb H,\gamma_2}(s) \] 
then \eqref{eq-prop2.4} ensures that for all $\tilde\gamma\in\Gamma$ we can take $\beta_{z_\psi}(\tilde\gamma)=\beta'_{z_{\psi'}}(\gamma\tilde\gamma\gamma^{-1})$. By Proposition \ref{lemma2.2}, the point $\mathcal P_{\mathbb H,\psi'}$ does not depend on the choice of a representative of $\boldsymbol\mu_{\mathbb H}$. Since $\mu^{\mathcal Y_{\rm rad}'}_{\mathbb H}$ is a representative of $\boldsymbol\mu_{\mathbb H}$ by \cite[Lemma 4.11]{LRV1}, it follows that 
\[ \mathcal P_{\mathbb H,\psi'}=t\cdot\beta'_{z_{\psi'}}(\gamma_{\psi'})=t\cdot\beta_{z_\psi}(\gamma_\psi)=\mathcal P_{\mathbb H,\psi} \] 
(observe that, by definition, $\gamma_{\psi'}=\gamma\gamma_\psi\gamma^{-1}$). This completes the proof. \end{proof}

\section{Formulas for Darmon points on $p$-adic tori} \label{sec-explicit}

\subsection{$p$-adic measures and Shapiro's lemma} \label{shapiro-subsec}

Define $\mathbb H_p:=\Z_p\otimes_\Z\mathbb H$. Let $\X:=(\Z_p^2)'$ be the set of primitive vectors in $\Y:=\Z_p^2$, i.e., those vectors in $\Y$ which are not divisible by $p$. Let $\tilde\D$ be the group of $\mathbb H_p$-valued measures on $\Y$ and write $\D$ for the subgroup consisting of those measures which are supported on $\X$. The group $\Sigma:=\GL_2(\Q_p)\cap\M_2(\Z_p)$ acts 
on $\tilde\D$ on the left by the integration formula 
\[ \int_\Y\varphi(x,y)d(\gamma\cdot\nu)(x,y)=\int_\Y\varphi\bigl(\gamma\cdot(x,y)\bigr)d\mu(x,y), \] 
with 
\[ \gamma\cdot(x,y):=(ax+by,cx+dy) \] 
for all step functions $\varphi:\Y\rightarrow\mathbb H_p$, all measures $\nu\in\tilde\D$ and all matrices $\gamma=\smallmat abcd\in\Sigma$. This action also induces an action of $\Sigma$ on $\D$ (see \cite[Lemma 7.4]{LRV1}). 

Denote by 
\[ \pi:\X\longrightarrow\PP^1(\Q_p) \] 
the fibration defined by $(a,b)\mapsto a/b$. By \cite[Theorem 7.5]{LRV1}, the canonical map 
\[ \pi_*:H^1(\Gamma_0,\D)\longrightarrow H^1(\Gamma_0,\mathcal M)\] 
induced by $\pi$ is surjective. Choose any $\boldsymbol{\tilde\mu}_{\mathbb H}\in H^1(\Gamma_0,\D)$ such that 
\[ \pi_*(\boldsymbol{\tilde\mu}_{\mathbb H})=\res_{\Gamma_0}(\boldsymbol\mu_{\mathbb H}), \]
where $\res_{\Gamma_0}$ is the canonical restriction map in cohomology. As explained in \cite[\S 7.3]{LRV1}, we can choose a representative $\tilde\mu_{\mathbb H}$ of $\boldsymbol{\tilde\mu}_{\mathbb H}$ such that $\pi_*(\tilde\mu_{\mathbb H,\gamma})=\mu_{\mathbb H,\gamma}$ for all $\gamma\in\Gamma_0$ (here $\pi_*$ stands for the map induced by $\pi$ between spaces of measures). Namely, 
\[ \mu_{\mathbb H,\gamma}(U)=\tilde\mu_{\mathbb H,\gamma}\bigl(\pi^{-1}(U)\bigr) \] 
for all $\gamma\in\Gamma_0$ and all compact open subsets $U\subset\PP^1(\Q_p)$. From now on we fix such a $\tilde\mu_{\mathbb H}$. 

For any compact open subset $U\subset\Y$ let $\D_U$ denote the subset of the measures on $\Y$ which are supported on $U$. Thus, in particular, $\D_\X=\D$. Define 
\[ \X_\infty:=\Z_p^\times\times p\Z_p,\qquad\X_{\rm aff}:=\Z_p\times\Z_p^\times, \] 
so that $\X=\X_\infty\coprod\X_{\rm aff}$. 

Shapiro's lemma provides canonical isomorphisms 
\[ \mathcal S:H^1(\Gamma_0,\D)\overset\simeq\longrightarrow H^1(\Gamma_1,\D_{\X_\infty}) \]
and 
\[ \hat{\mathcal S}:H^1(\hat\Gamma_0,\D_{\omega_p\X})\overset\simeq\longrightarrow H^1(\Gamma_1,\D_{\omega_p\X_\infty}). \]
Now we may consider the Hecke algebra $\T_1$ associated with the pair $(\Gamma_1,S_1)$ where 
\[ S_1:=B^\times\cap\prod_\ell S_\ell \] 
and the product is taken over all prime numbers $\ell$ (the local factors $S_\ell$ are defined before equation \eqref{S}). See \cite[\S2.2]{LRV1} for details. By definition, the Hecke operator $U_p\in\T_1$ is given by $U_p=\Gamma_1 g_0\Gamma_1$ for an element $g_0\in R_1$ of norm $p$. We also have the Atkin--Lehner involution $W_p=\Gamma_1\omega_p\Gamma_1\in\T_1$. Since the Hecke algebra $\T_1$ acts on $H^1(\Gamma_1,\D_\infty)$, we may use $\mathcal S$ and $\hat{\mathcal S}$ to define Hecke operators 
\[ \U_p:=\mathcal S^{-1}U_p\mathcal S:H^1(\Gamma_0,\D)\longrightarrow H^1(\Gamma_0,\D) \] 
and 
\[ \W_p:=\hat{\mathcal S}^{-1}W_p\mathcal S:H^1(\Gamma_0,\D)\longrightarrow H^1(\hat\Gamma_0,\D_{\omega_p\X}). \]
Define 
\[ \boldsymbol{\hat{\tilde\mu}}_{\mathbb H}:=\W_p\U_p(\boldsymbol{\tilde\mu}_{\mathbb H})\in H^1(\hat\Gamma_0,\D_{\omega_p\X}). \]
Choose a representative $\hat{\tilde\mu}_{\mathbb H}$ of $\boldsymbol{\hat{\tilde\mu}}_{\mathbb H}$. By \cite[eq. (45)]{LRV1}, we know that there are elements $m_1\in\D_{\X_\infty}$ and $m_2\in\D_{\omega_p\X_{\rm aff}}$ such that 
\[ \hat{\tilde\mu}_{\mathbb H,\gamma}=U_p^2\tilde\mu_\gamma+\gamma m_1-m_1\qquad\text{on $\X_\infty$} \] 
and 
\[ \hat{\tilde\mu}_{\mathbb H,\gamma}=p\tilde\mu_\gamma+\gamma m_2-m_2\qquad\text{on $p\X_{\rm aff}$} \]
for all $\gamma\in\Gamma_1$.

\subsection{An integral formula for Darmon points on tori} \label{integral-darmon-subsec}

We need to recall an auxiliary result from \cite{LV-darmon}. Let 
\[ G=G_1*_{G_0}G_2 \] 
be the amalgamated product of the groups $G_1$ and $G_2$ along the subgroup $G_0$. Let $A$ be an abelian group and consider the portion of the Mayer--Vietoris exact sequence given by 
\[ H^1(G_1,A)\oplus H^1(G_2,A)\longrightarrow H^1(G_0,A)\overset\Delta\longrightarrow H^2(G,A)\longrightarrow H^2(G_1,A)\oplus H^2(G_2,A). \] 
To simplify notations, in the next lemma we use the same symbols to denote cochains and cocycles representing them.
\begin{lemma} \label{lemma} 
Let $G$ be as above and assume that 
\begin{itemize} 
\item $c=\Delta(\rho)$ for some $\rho\in H^1(G_0,A)$, so there are $1$-cochains $\theta_1\in C^1(G_1,A)$ and $\theta_2\in C^1(G_2,A)$ such that $c_{|G_1}=\delta(\theta_1)$ and $c_{|G_2}=\delta(\theta_2)$, where $\delta$ is the connecting map and $\rho=\theta_{1|G_0}-\theta_{2|G_0}$;
\item $c$ is trivial in $H^2(G,A)$, so $c=\delta(b)$ for some $1$-cochain $b\in C^1(G,A)$.
\end{itemize}
Then $b_{|G_1}=\theta_1+\varphi_{1|G_1}$ and $b_{|G_2}=\theta_2+\varphi_{2|G_2}$ where $\varphi_1,\varphi_2:G\rightarrow A$ are group homomorphisms.
\end{lemma}

\begin{proof} An application of the universal property of the amalgamated product. For details, see \cite{LV-darmon}. \end{proof}  

From here on fix an algebraic closure $\bar\Q_p$ of $\Q_p$. Let $m\geq1$ be the rank of $\mathbb H$ and fix an isomorphism $\mathbb H\simeq\Z^m$, which allows us to identify $T_{\mathbb H}(\bar\Q_p)$ with $(\bar\Q_p^\times)^m$ and $\mathbb H_p$ with $\Z_p^m$. Let
\[ \Psi:T_{\mathbb H}(\bar\Q_p)/L_{\mathbb H}\longrightarrow\bar\Q_p^m \]
be a continuous homomorphism which is $K_p^m$-valued on $T_{\mathbb H}(K_p)/L_{\mathbb H}$. By pre-composing it with the projection onto the quotient (respectively, with the diagonal embedding and the projection onto the quotient) we will also view $\Psi$ as defined on $T_{\mathbb H}(\bar\Q_p)$ (respectively, on $\bar\Q_p^\times$). 

The first main result of this article is  

\begin{theorem} \label{thm2.11}
$\displaystyle{\Psi(\mathcal P_{\mathbb H,\psi})=-t\cdot\int_\X\Psi(x-z_\psi y)d{\tilde\mu}_{\mathbb H,\gamma_\psi}(x,y)}$. 
\end{theorem} 

\begin{proof} To begin with, observe that $x-z_\psi y\in(\mathcal O_K\otimes\Z_p)^\times$ when $(x,y)\in\X$. Define the $1$-cochain $\rho$ on $\Gamma_0$ with values in $T_{\mathbb H}(K_p)$ by the formula
\[ \rho_\gamma:=-\int_\X\Psi(x-z_\psi y)d\tilde\mu_{\mathbb H,\gamma}. \] 
Moreover, define the $1$-cochain $\hat\rho$ on $\hat\Gamma_0$ with values in $T_{\mathbb H}(K_p)$ by 
\[ \hat\rho_\gamma:=-\int_\X\Psi(x-z_\psi y)d\hat{\tilde\mu}_{\mathbb H,\gamma}+\int_{\X_\infty}\Psi(x-z_\psi y)d(\gamma m_1-m_1)+\int_{p\X_{\rm aff}}\Psi(x-z_\psi y)d(\gamma m_2-m_2). \]
Following \emph{verbatim} the proofs of \cite[Propositions 7.10 and 7.13]{LRV1} after replacing the branch $\log_p$ of the $p$-adic logarithm with the map $\Psi$, we see that $\rho$ and $\hat\rho$ split the restrictions of $\Psi(d)$ to $\Gamma_0$ and $\hat\Gamma_0$, respectively. Thus $\Psi(\boldsymbol{d}_z)=\Delta(\boldsymbol\rho-\boldsymbol{\hat\rho})$ where $\Delta$ is the connecting map in the Mayer--Vietoris exact sequence and $\boldsymbol{\rho}-\boldsymbol{\hat\rho}$ is the class in $H^1(\Gamma_1,T_{\mathbb H}(K_p))$ represented by $\rho-\hat\rho$. We already know that the $2$-cocycle $\Psi\bigl(\bar d_z\bigr)$ is split by $\Psi(\beta_z)$. Hence
\[ \delta\big(\Psi(\beta_z)\big)=\Delta(\rho-\hat\rho) \] 
where $\delta$ is the connecting map on cochains. Since $\Gamma=\Gamma_0*_{\Gamma_1}\hat\Gamma_0$ and $t$ annihilates $\Gamma^{\rm ab}$, the result follows from Lemma \ref{lemma}. \end{proof}

\section{Formulas for Darmon points on abelian varieties}

\subsection{Modular abelian varieties and isogenous tori} \label{modular-subsec}

From now on let $A_{/\Q}$ be a $d$-dimensional modular abelian variety of conductor $N:=MDp$. Therefore $A$ is attached to a (normalized) newform $f$ of weight $2$ and level $N$ via the Shimura construction. Namely, let $f(q)=\sum_{n\geq1}a_n(f)q^n$ be the $q$-expansion of $f$, let $F$ stand for the number field (of degree $d$) generated over $\Q$ by the Fourier coefficients $a_n(f)$ of $f$ and denote $\cO_F$ the ring of integers of $F$. Consider the homomorphism of algebras
\[ \lambda_f:\T\longrightarrow\cO_F \]
such that $\lambda_f(T_\ell)=a_\ell(f)$ whenever $\ell$ is a prime not dividing $N$. Then, as explained in \cite[\S 3.4.4]{zh}, $A$ is isogenous over $\Q$ to the maximal abelian subvariety of $J_1$ killed by $\ker(\lambda_f)$. Since $p|\!|N$, it is known that $A$ has purely multiplicative reduction at $p$, i.e., the identity component of the special fibre of the Néron model of $A$ over $\Z_p$ is a torus over $\F_p$ (cf. \cite[\S 14]{MTT}). Equivalently, the reduction of $A$ at $p$ is semistable of toric dimension $d$. Since $J_1^{\text{$p$-new}}$ is the maximal toric quotient of $J_1$ at $p$, it follows that $A$ is a quotient of $J_1^{\text{$p$-new}}$.

Set $\mathbb H_A:=H_1(A(\C),\Z)$, choose a sign $\epsilon\in\{\pm\}$ and write $\mathbb H_A^\epsilon$ for the $\epsilon$-eigenspace for complex conjugation acting on $\mathbb H_A$. Thus $\mathbb H_A^\epsilon$ is a free abelian group of rank $d$; fix an isomorphism $\mathbb H_A^\epsilon\simeq\Z^d$. Since $A$ is a quotient of $J_1^{\text{$p$-new}}$, it follows that $\mathbb H_A^\epsilon$ is a quotient of $H$.

Set $T_A:=T_{\mathbb H^\epsilon_A}$ and recall the Hecke-stable lattice $L_A:=L_{\mathbb H^\epsilon_A}$ inside $T_A(\Q_p)\simeq(\Q_p^\times)^d$. 

\begin{theorem} \label{homothetic-thm}
The $p$-adic torus $T_A/L_A$ and the abelian variety $A$ are isogenous over $K_p$.
\end{theorem}

\begin{proof} Multiplicity one ensures that $A$ is, up to isogeny, the unique quotient of $J_1^{\text{$p$-new}}$ on which the action of the Hecke operators $T_\ell$ for primes $\ell\nmid N$ factors through $\lambda_f$. Similarly, $\bar\Q_p^\times\otimes\mathbb H_A^\epsilon\simeq(\bar\Q_p^\times)^d$ is the unique quotient of $\bar\Q_p^\times\otimes H$ on which the action of these operators factors through $\lambda_f$ and complex conjugation acts as $\epsilon1$. The result follows from \cite[Theorem 1.1]{LRV1}. \end{proof}

\begin{remark} \label{definition-isogeny-rem}
The proof of \cite[Theorem 1.1]{LRV1} shows that the isogeny whose existence is ensured by Theorem \ref{homothetic-thm} can be chosen to be equivariant with respect to the actions of the relevant local Galois group (see \cite[\S 7.1]{LRV1}).
\end{remark}

\subsection{Darmon points on $A$} \label{darmon-A-subsec}

By Theorem \ref{homothetic-thm} and Remark \ref{definition-isogeny-rem}, we can fix a Galois-equivariant isogeny
\begin{equation} \label{varphi-A-eq}
\varphi_A:T_A(\bar\Q_p)/L_A\longrightarrow A(\bar\Q_p) 
\end{equation}
over $K_p$. For every $\psi\in\Emb(\cO_c,R_0)$ define
\[ P^\epsilon_{A,\psi}:=\varphi_A\bigl(\mathcal P_{\mathbb H^\epsilon_A,\psi}\bigr)\in A(K_p). \]
These are the Darmon points of conductor $c$ on $A$. Since the sign $\epsilon\in\{\pm\}$ has been chosen once and for all, from here only we shall simply write $P_{A,\psi}=P^\epsilon_{A,\psi}$. By Proposition \ref{conjugacy-prop}, $P_{A,\psi}=P_{A,\psi'}$ whenever $\psi$ and $\psi'$ are $\Gamma_0$-conjugate. In the next subsection we will give an integral expression for a suitable logarithm of the Darmon points $P_{A,\psi}$ on $A$.

\subsection{An integral formula for Darmon points on $A$}

Our goal is to describe a formula for Darmon points on $A$ which is a consequence of Theorem \ref{thm2.11}. In doing this, we replace the homomorphism $\Psi$ appearing in \S \ref{integral-darmon-subsec} with a map of particular arithmetic significance.

Let $\C_p$ be the completion of $\bar\Q_p$ and let $G$ be a finite-dimensional commutative Lie group over $\C_p$ (see \cite[Ch. III, \S 1]{Bou}). The Lie algebra $\text{Lie}(G)$ of $G$ is the tangent space of $G$ at the identity and is a $\C_p$-vector space of dimension $\text{dim}(G)$. Let $G_f$ be the smallest open subgroup of $G$ such that the quotient $G/G_f$ is torsion-free. As explained in \cite[Ch. III, \S 7.6]{Bou}, there is a canonical analytic homomorphism
\[ \log:G_f\longrightarrow\text{Lie}(G). \]
The map $\log$ is a local diffeomorphism and its kernel is the torsion subgroup of $G_f$.

In the special case where $G=A(\C_p)$, if $U$ is an open subgroup of $A(\C_p)$ then $A(\C_p)/U$ is torsion (\cite[Theorem 4.1]{Col}), hence $G_f=G$ and we get a logarithm map
\[ \log_A:A(\C_p)\longrightarrow\text{Lie}(A(\C_p))\simeq\C_p^d \]
on the $p$-adic points of $A$ whose kernel is the torsion subgroup of $A(\C_p)$. Fix an isomorphism $\text{Lie}(A(\C_p))\simeq\C_p^d$. As remarked in \cite[\S 1]{za}, $\log_A(A(H))\subset H^d$ for every complete subfield $H$ of $\C_p$, so that, in particular, $\log_A(P)\in K_p^d$ for all $P\in A(K_p)$.

\begin{remark}
It is known that $\log_A$ essentially coincides with the logarithm map introduced by Bloch and Kato in their fundamental paper \cite{BK} on the Tamagawa number conjecture for motives (cf. \cite[Example 3.11]{BK}).
\end{remark}

Composing $\log_A$ with the map $\varphi_A$ introduced in \eqref{varphi-A-eq} gives an analytic homomorphism
\[ \Psi_A:=\log_A\circ\,\varphi_A:T_A(\bar\Q_p)/L_A\longrightarrow\bar\Q_p^d \]
which is $K_p^d$-valued on $T_A(K_p)/L_A$.
 
Now let $f$ be an endomorphism of $A_{/\C_p}$ and let $f_*:\text{Lie}(A(\C_p))\rightarrow\text{Lie}(A(\C_p))$ be the linear map induced by $f$. By \cite[\S 1]{za}, the maps $\log_A$ and $f$ commute, in the sense that
\begin{equation} \label{comm-log-eq}
\log_A\circ\,f=f_*\circ\log_A.
\end{equation}
With notation as in \S \ref{shapiro-subsec}, for simplicity set $\tilde\mu_A:=\tilde\mu_{\mathbb H_A^\epsilon}$. Let $\psi\in\Emb(\cO_c,R_0)$ and recall the Darmon point $P_{A,\psi}\in A(K_p)$ defined in \S \ref{darmon-A-subsec}. The second main result of this paper is 

\begin{theorem} \label{thm2}
$\displaystyle{\log_A\bigl(f(P_{A,\psi})\bigr)=-t\cdot f_*\left(\int_\X\Psi_A(x-z_\psi y)d{\tilde\mu}_{A,\gamma_\psi}(x,y)\right)}$. 
\end{theorem}

\begin{proof} Immediate from Theorem \ref{thm2.11} upon taking $\mathbb H=\mathbb H_A^\epsilon$, $\Psi=\Psi_A$ and using \eqref{comm-log-eq}. \end{proof}

By choosing $f=\text{id}_A$ in Theorem \ref{thm2}, we deduce

\begin{corollary} \label{thm2-coro2}
$\displaystyle{\log_A(P_{A,\psi})=-t\cdot\int_\X\Psi_A(x-z_\psi y)d{\tilde\mu}_{A,\gamma_\psi}(x,y)}$.
\end{corollary}

While the validity of Theorem \ref{thm2} is unconditional, now we describe a refinement that depends on a technical assumption. As shown in \cite[Ch. III, \S 7, Proposition 10]{Bou}, there is an open neighbourhood $V$ of the identity in $A(\C_p)$ such that ${\log_A|}_V$ is the inverse of the restriction of the exponential map
\[ \exp_A:\text{Lie}(A(\C_p))\simeq\C_p^d\longrightarrow A(\C_p) \]
to a suitable open subgroup of $\text{Lie}(A(\C_p))$. 

\begin{corollary} \label{thm2-coro}
With notation as before, assume that $f(P_{A,\psi})\in V$. Then
\[ f(P_{A,\psi})=\exp_A\!\Bigg(\!\!-t\cdot f_*\left(\int_\X\Psi_A(x-z_\psi y)d{\tilde\mu}_{A,\gamma_\psi}(x,y)\right)\!\!\Bigg). \]
\end{corollary}

\begin{proof} Since $f(P_{A,\psi})\in V$, the claim follows by applying $\exp_A$ to both sides of the equality in Theorem \ref{thm2}. \end{proof}

In other words, Corollary \ref{thm2-coro} says that if the point $f(P_{A,\psi})$ is ``sufficiently close'' to the identity of $A$ then we can exhibit a formula not only for the logarithm of $f(P_{A,\psi})$ but also for $f(P_{A,\psi})$ itself.

\appendix

\section{The $p$-adic uniformization of abelian varieties} \label{appendix}

We review the $p$-adic uniformization theory of abelian varieties $A$ over $\Q$ with purely split multiplicative reduction at $p$ and, under this assumption, we describe an explicit construction of the isogeny $\varphi_A$ appearing in \S \ref{darmon-A-subsec}.

\subsection{Review of Tate--Morikawa--Mumford theory} \label{tate-subsec}

Let $A_{/\Q}$ be an abelian variety of dimension $d$ with purely \emph{split} multiplicative reduction at $p$. This means that the identity component of the special fibre of the Néron model of $A$ over $\Z_p$ is a \emph{split} torus over $\F_p$. Such an abelian variety is a higher-dimensional analogue of Tate's $p$-adic elliptic curve (see, e.g., \cite[\S A.1.1]{se2}). Thanks to results due (among others) to Tate, Morikawa and Mumford (\cite{mc}, \cite{mo1}, \cite{mo2}, \cite{mu2}), which now we briefly describe, this condition guarantees that $A$ admits an analytic uniformization locally at $p$. Details can be found in \cite[Section III]{ri}, while the rigid analytic point of view is nicely exposed in \cite[Section 4]{pa}.

Set $G_{\Q_p}:=\Gal(\bar\Q_p/\Q_p)$. There are free abelian groups $\mathcal M$ and $\mathcal N$ of rank $d$ and an admissible (in the sense of \cite[Section III, \S 2]{ri}) homomorphism $\alpha:\mathcal M\rightarrow\Hom(\mathcal N,\Q_p^\times)$ that fit into a short exact sequence
\[ 0\longrightarrow\mathcal M\overset\alpha\longrightarrow\Hom(\mathcal N,\bar\Q_p^\times)\overset\theta\longrightarrow A(\bar\Q_p)\longrightarrow0 \]
of $G_{\Q_p}$-modules (with $\mathcal M$ and $\mathcal N$ being regarded as trivial Galois modules). In other words, there is a Galois-equivariant analytic uniformization
\[ \theta:\Hom(\mathcal N,\bar\Q_p^\times)/\mathcal M\overset\simeq\longrightarrow A(\bar\Q_p) \]
which expresses the geometric points of $A_{/\Q_p}$ as a quotient of a $d$-dimensional $p$-adic torus by a sublattice of full rank. Reversing the roles of $\mathcal M$ and $\mathcal N$, one obtains an analogous parametrization for the dual abelian variety of $A$, but we will not need this fact. 

Choose $\Z$-bases $\{x_1,\dots,x_d\}$ and $\{y_1,\dots,y_d\}$ of $\mathcal M$ and $\mathcal N$, respectively, and for every $j=1,\dots,d$ define
\[ \q_j:=\bigl(\alpha(x_j)(y_1),\dots,\alpha(x_j)(y_d)\bigr)\in(\Q_p^\times)^d. \]
It follows that $\theta$ induces an analytic isomorphism
\begin{equation} \label{analytic-iso-eq}
\Phi_{\rm Tate}:(\bar\Q_p^\times)^d\big/\big\langle\q_1,\dots,\q_d\big\rangle\overset\simeq\longrightarrow A(\bar\Q_p) 
\end{equation}
of $G_{\Q_p}$-modules. Since we are assuming that the reduction of $A$ at $p$ is split, the map $\Phi_{\rm Tate}$ is defined over $\Q_p$. The vectors $\q_1,\dots,\q_d$ are the $d$-dimensional analogue of the period $q$ which appears in the theory of Tate's elliptic curves (see, e.g., \cite[\S A.1.1]{se2}), and they will henceforth be referred to as \emph{Tate periods} for $A$ at $p$.

\begin{remark}
If the reduction of $A$ at $p$ is purely multiplicative but not necessarily split then a result of Mumford and Raynaud shows the existence of an unramified extension $H$ of $\Q_p$ such that $A$ has a $p$-adic uniformization \eqref{analytic-iso-eq} defined over $H$ (cf. \cite[Theorem 3.2.2]{ri}).
\end{remark}

\subsection{Modularity and homothetic lattices}

Now assume that the $d$-dimensional abelian variety $A_{/\Q}$ has conductor $MDp$ and is modular. As remarked in \S \ref{modular-subsec}, it is known that $A$ has purely multiplicative reduction at $p$, and we suppose, in addition, that the reduction is \emph{split}, so that the results reviewed in \S \ref{tate-subsec} apply. 

Recall the $d$-dimensional $p$-adic torus $T_A$ and the Hecke-stable lattice $L_A$ inside $T_A(\Q_p)$ defined in \S \ref{modular-subsec}. As in \cite[Definition 29]{Gr}, we say that two lattices $\Lambda_1$ and $\Lambda_2$ in $(K_p^\times)^d$ are \emph{homothetic} if $\Lambda_1\cap\Lambda_2$ has finite index both in $\Lambda_1$ and in $\Lambda_2$. 

The following result, asserting that $L_A$ is homothetic to the lattice of Tate periods for $A$ at $p$, extends \cite[Theorem 7.16]{LRV1} from elliptic curves to abelian varieties of higher dimension.

\begin{theorem} \label{homothetic2-thm}
The lattices $L_A$ and $\langle\q_1,\dots,\q_d\rangle$ are homothetic in $(K_p^\times)^d$.
\end{theorem}

\begin{proof} Arguing as in the proof of Theorem \ref{homothetic-thm}, this is a consequence of multiplicity one and the uniformization results of \cite{LRV1}. \end{proof}

\subsection{The isogeny $\varphi_A$}

By Theorem \ref{homothetic2-thm}, $L_A\cap\langle\q_1,\dots,\q_d\rangle$ has finite index in $L_A$; set 
\[ n:=\bigl[L_A:L_A\cap\langle\q_1,\dots,\q_d\rangle\bigr]. \]
Raising to the $n$-th power gives a Galois-equivariant isogeny of $p$-adic tori
\begin{equation} \label{p-adic-isogeny-eq}
T_A(\bar\Q_p)/L_A\longrightarrow(\bar\Q_p^\times)^d\big/\big\langle\q_1,\dots,\q_d\big\rangle
\end{equation}
defined over $\Q_p$, and then composing \eqref{p-adic-isogeny-eq} with $\Phi_{\rm Tate}$ yields a Galois-equivariant isogeny
\[ \varphi_A:T_A(\bar\Q_p)/L_A\longrightarrow A(\bar\Q_p). \]
This gives an explicit construction of the isogeny introduced in \eqref{varphi-A-eq}.

\end{document}